\documentclass[12pt]{amsart}

\usepackage{a4}   
\usepackage{amsmath,amssymb,amsthm,mathrsfs,textcomp}   
\usepackage{amsfonts}
\usepackage{mathtools}   
\usepackage{wasysym}   
\usepackage{csquotes}   
\usepackage{marvosym}

\newtheorem{thm}{Theorem}
\newtheorem{lem}{Lemma}
\newtheorem{prop}{Proposition}
\newtheorem{cor}{Corollary}
\theoremstyle{remark}

\newtheorem*{rmk_plain}{Remark}

\newcommand{\N}{\mathbb{N}}

\newcommand{\R}{\mathbb{R}}
\newcommand{\C}{\mathbb{C}}
\newcommand{\D}{\mathbb{D}}

\newcommand{\T}{\mathbb{T}}

\newcommand{\coloneqq}{\mathrel{\mathop:}=}

\DeclareMathOperator{\RE}{Re}
\DeclareMathOperator{\SUPP}{supp}
\DeclareMathOperator{\DIST}{dist}

\makeatletter
\def\moverlay{\mathpalette\mov@rlay}
\def\mov@rlay#1#2{\leavevmode\vtop{%
   \baselineskip\z@skip \lineskiplimit-\maxdimen
   \ialign{\hfil$\m@th#1##$\hfil\cr#2\crcr}}}
\newcommand{\charfusion}[3][\mathord]{
    #1{\ifx#1\mathop\vphantom{#2}\fi
        \mathpalette\mov@rlay{#2\cr#3}
      }
    \ifx#1\mathop\expandafter\displaylimits\fi}
\makeatother

\makeatletter
\let\save@mathaccent\mathaccent
\newcommand*\if@single[3]{%
  \setbox0\hbox{${\mathaccent"0362{#1}}^H$}%
  \setbox2\hbox{${\mathaccent"0362{\kern0pt#1}}^H$}%
  \ifdim\ht0=\ht2 #3\else #2\fi
  }
\newcommand*\rel@kern[1]{\kern#1\dimexpr\macc@kerna}
\newcommand*\widebar[1]{\@ifnextchar^{{\wide@bar{#1}{0}}}{\wide@bar{#1}{1}}}
\newcommand*\wide@bar[2]{\if@single{#1}{\wide@bar@{#1}{#2}{1}}{\wide@bar@{#1}{#2}{2}}}
\newcommand*\wide@bar@[3]{%
  \begingroup
  \def\mathaccent##1##2{%
    \let\mathaccent\save@mathaccent
    \if#32 \let\macc@nucleus\first@char \fi
    \setbox\z@\hbox{$\macc@style{\macc@nucleus}_{}$}%
    \setbox\tw@\hbox{$\macc@style{\macc@nucleus}{}_{}$}%
    \dimen@\wd\tw@
    \advance\dimen@-\wd\z@
    \divide\dimen@ 3
    \@tempdima\wd\tw@
    \advance\@tempdima-\scriptspace
    \divide\@tempdima 10
    \advance\dimen@-\@tempdima
    \ifdim\dimen@>\z@ \dimen@0pt\fi
    \rel@kern{0.6}\kern-\dimen@
    \if#31
      \overline{\rel@kern{-0.6}\kern\dimen@\macc@nucleus\rel@kern{0.4}\kern\dimen@}%
      \advance\dimen@0.4\dimexpr\macc@kerna
      \let\final@kern#2%
      \ifdim\dimen@<\z@ \let\final@kern1\fi
      \if\final@kern1 \kern-\dimen@\fi
    \else
      \overline{\rel@kern{-0.6}\kern\dimen@#1}%
    \fi
  }%
  \macc@depth\@ne
  \let\math@bgroup\@empty \let\math@egroup\macc@set@skewchar
  \mathsurround\z@ \frozen@everymath{\mathgroup\macc@group\relax}%
  \macc@set@skewchar\relax
  \let\mathaccentV\macc@nested@a
  \if#31
    \macc@nested@a\relax111{#1}%
  \else
    \def\gobble@till@marker##1\endmarker{}%
    \futurelet\first@char\gobble@till@marker#1\endmarker
    \ifcat\noexpand\first@char A\else
      \def\first@char{}%
    \fi
    \macc@nested@a\relax111{\first@char}%
  \fi
  \endgroup
}
\makeatother

\makeatletter
\renewcommand*\env@matrix[1][*\c@MaxMatrixCols c]{%
  \hskip -\arraycolsep
  \let\@ifnextchar\new@ifnextchar
  \array{#1}}
\makeatother

\def\clap#1{\hbox to 0pt{\hss#1\hss}}

\begin{document}
\title[Hyperconvexity of pseudoconvex domains]{An answer to \enquote{Hyperconvexity of pseudoconvex domains with continuous boundaries} by X. Wang and C. Yang}
\author{Tobias Harz}
\address{Bergische Universität Wuppertal, Wuppertal, Germany}
\email{harz@math.uni-wuppertal.de}
\date{\today}

\subjclass[2020]{Primary 32U10; Secondary 32T35, 32U05, 32U35}
\keywords{Hyperconvex domains, continuous boundaries, Poisson envelopes, deformation of analytic discs}
\begin{abstract}
We give a variant of the recent proof concerning hyperconvexity of pseudoconvex domains with continuous boundary.
\end{abstract}
\maketitle

\begin{thm}
Every bounded pseudoconvex domain $\Omega \subset \C^n$ with continuous boundary is hyperconvex.
\end{thm}
\begin{proof}
Consider coordinates $(w,z) \in \C \times \C^{n-1}$, $w = u+iv$, and write
$$
 G \coloneqq \{(v,z) \in \R \times \C^{n-1}: v^2 + |z|^2 < 1\}, \quad \hat{G} \coloneqq \{(v,z) \in \R \times \C^{n-1}: v^2 + |z|^2 < 2\}.
$$
Recall that hyperconvexity is a local property of the boundary, see Proposition 1.1 in \cite{KerzmanRosay81}, and that a domain $\Omega \subset\subset \C^n$ is hyperconvex if and only if it admits a (not necessarily continuous) negative plurisubharmonic (psh) exhaustion function, see Theorem 1.6 in \cite{Blocki96}. It thus suffices to consider domains of the form
\begin{align*}
 \Omega \coloneqq \{ (w,z) \in \R_u \times G: 0 < u < F(v,z) \}, \quad  \hat{\Omega} &\coloneqq \{ (w,z) \in \R_u \times \hat{G}: u < F(v,z)\},
\end{align*}
where $F \colon \hat{G} \to [1,\infty)$ is a continuous function such that $\hat\Omega$ is pseudoconvex, and to show that there exists a psh function 
$\varphi \colon \Omega \to (-\infty,0)$ with $\lim_{q \to p} \varphi(q) = 0 \text{ for every } p \in \Gamma \coloneqq \{(w,z) \in \R_u \times G : u = F(v,z)\}$ (since then $(w,z) \mapsto \max\{\varphi(w,z), v^2+|z|^2-1,-u\}$ is a negative psh exhaustion function on $\Omega$).

Define $r \colon \Omega \to \R$ as $r(w,z) \coloneqq u-F(v,z)$, and consider its Poisson envelope $\hat{r} \colon \Omega \to \R$ given by
$$
\hat{r}(q) \coloneqq \frac{1}{2\pi} \inf_{f \in \mathcal{F}(\Omega,q)} \int_\T (r \circ f) d\sigma,
$$
where $\mathcal{F}(\Omega,q) = \{f \in \mathcal{O}(\bar\D, \Omega) : f(0) = q\}$, and $\sigma$ denotes normalized arclength measure on $\T \coloneqq \{\zeta \in \C : |\zeta| = 1\}$. Obviously, $\hat{r} < 0$ on $\Omega$.
Moreover, it was shown by Poletsky \cite{Poletsky91} that $\hat{r}$ is psh (and, in fact, the largest psh minorant of $r$). 
Hence, if one can show $\lim_{q \to p} \hat{r}(q) = 0$ for every $p \in \Gamma$, then $\varphi \coloneqq \hat{r}$ is a function as desired. 

Assume, in order to get a contradiction, that there exists $p \in \Gamma$ such that $\liminf_{q \to p} \hat{r}(q) < 0$. 
Then there exist $\varepsilon \in (0,1)$, a sequence $\{q_j\} \subset \Omega$ converging to $p$, and analytic discs $f_j \in \mathcal{F}(\Omega,q_j)$ such that
\begin{align} \label{E:contradiction}
\frac{1}{2\pi} \int_\T (r \circ f_j) d\sigma < -\varepsilon \;\; \text{for every } j \in \N.
\end{align}
While the centers $f_j(0)$ converge to the boundary point $p$, the existence of the sequence $(f_j)$ is not a contradiction to the continuity principle, 
since the boundaries $f_j(b\D)$ do not have to stay in a compact subset of $\Omega$. However, in view of the definition of $r$, condition \eqref{E:contradiction} 
states that a uniformly large part of the boundary of each disc stays uniformly far away from $\Gamma$ in the coordinate direction $u$. 
The goal is to exploit this property in order to perturb the discs $(f_j)$ in such a way that a contradiction to the continuity principle is reached.

For every $j \in \N$, consider the open set $E_j \coloneqq \{\zeta \in \T : (r \circ f_j)(\zeta) < -\frac{\varepsilon}{2}\}$, and note that in view of \eqref{E:contradiction} and boundedness of $\Omega$ there exists $\delta > 0$ such that $\sigma(E_j) > \delta$ for all $j \in \N$. Hence, for some fixed sufficiently small $m>0$, we can choose 
$$
u_j \colon \T \to [0,\textstyle\frac{\varepsilon}{2}\displaystyle] \text{ smooth}, \quad \SUPP u_j \subset E_j, \quad \frac{1}{2\pi}\int_\T u_j d\sigma = m \text{ for all } j \in \N.
$$
There exist unique functions
$$
h_j = \hat{u}_j + i\hat{v}_j \in \mathcal{O}(\D) \cap \mathcal{C}(\bar\D) \quad\text{such that}\quad \hat{u}_j|_\T = u_j, \quad h_j(0) = m.
$$
The naive idea is to consider the perturbations $\tilde{f}_j(\zeta) \coloneqq f_j(\zeta) + (h_j(\zeta),0')$, which satisfy $\tilde{f}_j(0) = q_j + (m,0')$, with $m>0$ independent of $j \in \N$. However, while the distance of $\tilde{f_j}(\T)$ to $\Gamma$ in the $u$-direction is controlled, namely $f_j(\zeta) + (u_j(\zeta),0') \in \Omega$ for $\zeta \in \T$, the problem is that the Hilbert transform $v_j \coloneqq \hat{v}_j|_\T$ of $u_j$ need not be uniformly bounded, and in particular $\tilde{f}_j(\T)$ need not be contained in $\hat\Omega$. It is precisely here that Wang and Yang \cite{WangYang26} introduce a crucial observation, by means of a suitably chosen holomorphic transformation $\Phi$.

Namely, the function $\Phi(\lambda) \coloneqq \frac{1}{3}(\cosh \lambda - 1)$, $\lambda = x+iy \in \C$, has the following properties\footnote{Note that $\cosh \lambda-1 = -2J(e^{-\lambda})$, where $J(\lambda) = \frac{2-\lambda-\lambda^{-1}}{4}$ is the Joukowski function that appears in \cite{WangYang26}. To check the listed properties, compute with $\lambda=x+iy$: \begin{enumerate} \item[(i)] $|\Phi(\lambda)| \le \frac{1}{3}(|\frac{e^\lambda+e^{-\lambda}}{2}| + 1) \le \frac{1}{3}(\frac{e^x+e^{-x}}{2} + 1) \le \frac{1}{3}(\cosh 1 + 1) < 1$ for $x \in [0,1]$. \item[(ii)] $\Phi(iy) = \frac{1}{3}\Big(\frac{e^{iy}+e^{-iy}}{2}-1\Big) = \frac{1}{3}(\cos y -1) = -\frac{2}{3}\sin^2(\frac{y}{2})$. \item[(iii)] $\Phi(x) = \frac{1}{6}(e^x + e^{-x}-2) = \frac{1}{6}(e^{x/2}-e^{-x/2})^2$. \end{enumerate}}:
\begin{enumerate}
 \item[(i)] $|\Phi| \le 1$ on $\{0 \le \RE \lambda \le 1\}$,
 \item[(ii)] $\Phi(iy) \in (-\infty,0]$, $y \in \R$,
 \item[(iii)] $\Phi(x) \in (0,\infty)$, $x \in \R \setminus\{0\}$.
\end{enumerate}
Since $F$ is continuous, we have $\{q \in \Omega: r(q)<-\frac{\varepsilon}{2}\} \subset\subset \hat\Omega$, and thus there exists $s \in (0,1)$ such that $\DIST(q,b\hat\Omega)>s$ on $\bigcup_{j\in\N}f_j(E_j)$.
Hence if
$$
\tilde{f}_{j,\tau}(\zeta) \coloneqq f_j(\zeta) + (\tau s(\Phi \circ h_j)(\zeta),0'), \quad j\in\N, \quad\tau\in[0,1],
$$
then for every $j \in \N$ and $\tau \in [0,1]$ we have
$$
 \tilde{f}_{j,0}(\bar\D) \subset \hat\Omega,\quad \tilde{f}_{j,\tau}(\T) \subset \hat\Omega, \quad \tilde{f}_{j,\tau}(0) = q_j + (\tau c,0'),
$$
where $c \coloneqq s\Phi(m) > 0$ is independent of $j$ (for the second statement, note that $u_j \equiv 0$ on $\T \setminus E_j$ and thus $(\Phi \circ h_j)(\zeta) \in (-\infty,0]$ for $\zeta \in \T \setminus E_j$ by (ii) -- this is precisely the key observation that allows us to control the Hilbert transform of $u_j$ at points $\zeta \in \T$ for which $f_j(\zeta)$ is close to $\Gamma$).
Since $q_j \to p$, we can fix $j \in \N$ so large that $r(q_j) > -c$, and thus $(\tilde{f}_{j,\tau})_{\tau \in [0,1]}$ is a $1$-parameter family of analytic discs violating the continuity principle for $\hat\Omega$. \Lightning
\end{proof}

\begin{rmk_plain}
Consider $\Omega \subset \hat\Omega$ as in the proof of the theorem. Suppose that we already know that $\hat\Omega$ is hyperconvex, i.e.\ there exists a psh exhaustion function $\sigma \colon \hat\Omega \to [-1,0)$. We claim that we can construct a continuous increasing convex function $\chi \colon [-1,0] \to (-\infty,0]$ such that 
\begin{align} \label{E:estimate}
 \chi(0) = 0 \quad\text{and}\quad \chi \circ \sigma \le r \quad\text{on } \Omega.
\end{align}
Since $\chi \circ \sigma$ is psh and $\hat{r}$ is the largest psh minorant of $r$, it follows that $\chi \circ \sigma \le \hat{r} < 0$ on $\Omega$, and since $\lim_{q \to p} (\chi \circ \sigma)(q) = 0$ for every $p \in \Gamma$, this implies that $\lim_{q \to p} \hat{r}(q) = 0$ for every $p \in \Gamma$. This proves that whenever $\hat\Omega$ is hyperconvex, $\hat{r}$ is a function as desired.

The author had long been aware that the function $\hat{r}$ was the natural function to investigate, but had been unable to prove $\lim_{q \to p} \hat{r}(q) = 0$ for $p \in \Gamma$ when $\hat\Omega$ is pseudoconvex and $F$ is continuous. The problem was precisely the lack of control over the Hilbert transforms $\{v_j\}$. In the author's view, the idea of controlling the Hilbert transform by postcomposing with the map $\Phi$ is the crucial ingredient that appears in \cite{WangYang26}.

To establish \eqref{E:estimate}, set $-M \coloneqq \inf_\Omega r > -\infty$, and define
$$
\omega(t) \coloneqq \begin{cases} \displaystyle\inf_{q \in \Omega \atop t \le \sigma(q)}r(q) & t \in [-1,0) \\ 0 & t=0 \end{cases}.
$$
Clearly, $\omega \colon [-1,0] \to [-M,0]$ is nondecreasing and satisfies $\omega \circ \sigma \le r$ on $\Omega$. Moreover, $\lim_{t \to 0^-} \omega(t) = 0$, since otherwise there would exist $\varepsilon > 0$ and a sequence $\{q_j\} \subset \Omega$ such that $\sigma(q_j) \to 0$ and $r(q_j) \le -\varepsilon$. But, by continuity of $F$, we have $K_\varepsilon	\coloneqq \{q \in \Omega : r(q) \le -\varepsilon\} \subset\subset \hat\Omega$, and thus $\sup_{K_\varepsilon}\sigma < 0$, contradicting $\sigma(q_j) \to 0$. Then let $\chi \colon [-1,0] \to (-\infty,0]$ be a continuous convex increasing minorant of $\omega$ such that $\chi(0) = 0$ (for example, take $\chi(t) \coloneqq \sup_{-1\le s < 0}\{\omega(s) +  \frac{M}{|s|}t\}$).
\end{rmk_plain}

\end{document}

The Poisson envelope $\hat{r}$ of $r(z,w) = u - F(v,z)$ is a natural candidate to construct a negative psh exhaustion function in the local setting described in the proof. 

Namely, as a proof of concept assume that the function $F$ in the above construction is $\mathcal{C}^2$-smooth. It follows easily from the classical result of Diederich-Forn{\ae}ss that there exists $\rho \in \mathcal{C}^2(\hat{Z})$, $Z \coloneqq \R_u \times \hat{G}$, such that $\hat\Omega = \{\rho<0\}$ and such that for some $\eta \in (0,1)$ the function $-(-\rho)^\eta$ is psh on $\Omega$. Since also $\hat\Omega = \{r<0\}$, where we view $r$ as a function $r \in \mathcal{C}^2(\hat{Z})$, it follows that $r=h\rho$ for a function $h \in \mathcal{C}^1(\hat{Z})$. In particular, for $C>0$ large enough, the function $\psi \coloneqq (Ch\rho)|_\Omega$ will be a psh function on $\Omega$ such that $\psi < r$ on $\Omega$. Since $\hat{r}$ is the largest psh minorant of $r$, we conclude that $\psi \le \hat{r} < 0$. But $\lim_{q \to p} \psi(q) \psi(q) = 0$ for every $p \in \Gamma$, which implies that $\lim_{q \to p}\hat{r}(q) = 0$ for $\Gamma$. This simple observation shows that, in the more special case of  it is already clear form the Diederich-Fornaess